\documentclass[12pt, reqno]{amsart}
\usepackage{amsmath, amstext, amsbsy, amssymb, amscd}

\newtheorem{theorem}{Theorem}[section]

\newtheorem{proposition}[theorem]{Proposition}
\newtheorem{corollary}[theorem]{Corollary}
\theoremstyle{definition}
\newtheorem{definition}[theorem]{Definition}

\theoremstyle{remark}
\newtheorem{remark}[theorem]{Remark}

\numberwithin{equation}{section}

{\vskip-\lastskip\medskip
  \noindent
  {\em #1.}\enspace
  }%
{\qed\par\medskip
  }

\begin{document}
\title
[Bertini Type Theorems]
{  Bertini Type Theorems}

\author[Jing  Zhang]{Jing  Zhang}
\address{ Department of Mathematics and Statistics, 
University at  Albany, SUNY, Albany, NY 12222,  USA}
\email{jzhang@albany.edu}

\begin{abstract}

Let $X$  be a smooth irreducible projective 
variety of dimension at least 2 over an
algebraically closed field of characteristic 0
in the projective space ${\mathbb{P}}^n$. 
 Bertini's Theorem states that 
a  general  hyperplane  $H$
intersects  $X$  with an irreducible 
smooth subvariety  of  $X$. However, the precise location of 
the smooth hyperplane section is 
not known. 
We show that   for any $q\leq n+1$ closed points  
in general position
and any degree $a>1$,
a general    hypersurface 
$H$  of  degree  $a$ passing through  these $q$ points intersects $X$ with  an irreducible  smooth  codimension 1  subvariety  
on  $X$. We also consider linear system of ample divisors and 
give precise location of smooth 
elements in the system. Similar result can be obtained for compact
complex manifolds with holomorphic maps into projective spaces. 
\end{abstract} 

\maketitle
\begin{flushleft}
\end{flushleft}
\date{}

\section{Introduction}

Bertini's two fundamental 
 theorems  concern the irreducibility and 
 smoothness of the general hyperplane section 
 of a smooth projective variety and a general member 
 of a linear system 
 of divisors. 
 The hyperplane version of  Bertini's theorems says that if $X$
is a smooth irreducible projective variety
of dimension at least 2
over an algebraically closed field
$k$ of characteristic 0, then  
a general hyperplane  $H$ intersects 
$X$ with a smooth irreducible subvariety 
of codimension 1 on $X$. 
But we do not know the exact location of the 
smooth hyperplane sections.

Let $F$  be an effective  divisor on $X$.  We say that $F$
is   a   fixed  component   of  linear   system 
$|D|$  of a divisor $D$  if   $E>F$   for all  $E\in  |D|$. 
$F$  
is the fixed part of a linear system if  every irreducible component of $F$ 
is a fixed component of the system and 
$F$  is maximal with respect to 
the order $\geq$. Every element $E$ in the system 
can be written in the form $E=E'+F$. We say that $E'$ is the 
variable part of $E$. A point $x\in X$  is a base point 
of the linear system   if $x$  is contained in the supports of
 variable parts of all divisors in the system. 
 The second Bertini Theorem is (\cite{U}, Theorem 4.21):  
If $\kappa(D, X)\geq 2$, then
the variable part of a general 
member of the complete linear system
$|D|$ is irreducible and smooth away from the singular locus 
of $X$ and the base locus of $|D|$. Here 
$$\kappa(D, X)=tr.deg_{\mathbb{C}}\oplus_{m\geq 0}
H^0(X, {\mathcal{O}}_X(mD))-1.$$

In \cite{F}, by using  the theory of
intersection  numbers  of  semipositive   line bundles, 
 Fujita sharpens the
above celebrated Bertini's theorem and
presents conditions on the base locus
 under which  the general member is also nonsingular 
 on the base locus itself. 
 In \cite{X}, Xu  applies deformation 
 of singularities  to
study the singularities of a generic element 
of a linear system and give  detailed information 
on the singular type of the base element.  In \cite{Z},  Zak
considers  that under what condition the hyperplane section 
of a normal variety is normal.  Diaz and Harbater consider the 
 singular locus of the
general member of a linear system and obtain better 
dimension estimate if the base locus is scheme-theoretically smooth.
They successfully apply their strong Bertini  theorem 
to complete intersection varieties.
Our results and methods
are    different from all these known results.
This work is   inspired by  Hartshorne's proof 
(\cite{H} Theorem 8.18, Chapter 2) and 
Kleiman's very interesting article
\cite{K1}.

In  this paper, we assume that the ground field $k$  is algebraically 
closed and of characteristic  0. 

\begin{definition} Let $S=\{P_0$, $P_1$,..., $P_{q-1}\}$ be $q$ points in ${\mathbb{P}}^n$.
We say that they are in general position if

(1) for $q<n+1$, the vectors defined by the homogeneous 
coordinates  of these  $q$  points  are linearly independent;

(2) for $q=n+1$, any $n$ points are linearly independent.

\end{definition}

Let  $L$  be the linear system 
of hypersurfaces  of  degree 
$a>1$  passing  through  these 
$q$  points   $P_0, P_1,...,  P_{q-1}$
in general position. 
Our main result  is 
that  a general  member of  $L$  is irreducible and smooth.

\begin{theorem}
If  $X$  is an irreducible   smooth projective  variety
of dimension at least  2 in ${\mathbb{P}}^n$, then  for any $q\leq n+1$
closed points 
 $P_0,$ $P_1,..., P_{q-1}$  on  $X$
 in a general position
and any degree $a>1$,
a general    hypersurface 
$H$  of  degree  $a$ passing through  these $q$  points 
intersects $X$  with 
  an irreducible  smooth  codimension 1  subvariety  
on  $X$.
\end{theorem}

In fact, if some points even all points do not lie on $X$,
Theorem 1.2  still  holds.

\begin{theorem}
If  $X$  is an irreducible   smooth projective  variety
of dimension at least  2 in ${\mathbb{P}}^n$, $D$  is an ample  divisor 
on  $X$, then there is an $n_0>0$ such that for all $m\geq n_0$ 
and any $q\leq n+1$ closed points $P_0, P_1,..., P_{q-1}$ on $X$, 
a general member of $|mD|_{q}$ is irreducible and  smooth, 
where $|mD|_{q}$
is the linear system of effective divisors in
$|mD|$ passing through 
these  $q$ points $P_s$, $s=0,1,..., (q-1)$.
\end{theorem}

The paper is organized as follows.
In Section 2 and 3, we will deal with hypersurface sections.
In Section 4 and 5, we will consider linear system of ample and big divisors.
In Section 6,  some  applications  in  compact  complex  manifolds
 will be discussed.  

\section{Hypersurface Sections Passing Through a Point}

Because of the lengthy  calculation,
we   first  show   the case when there is only one point 
to indicate  the idea. 
The general case will be proved in Section 3.

\begin{theorem}
If  $X$  is an irreducible   smooth projective  variety
of dimension at least  2, then  for any closed point  $P_0$  on  $X$
and any degree $a>1$,
a general    hypersurface 
$H$  of  degree  $a$ passing through  $P_0$ 
intersects $X$  with 
  an irreducible  smooth  codimension 1  subvariety  
on  $X$.
\end{theorem}

$Proof$. Let  $X$  be  a closed subset of  ${\mathbb{P}}^n$,
$n\geq  3$.
We may assume     that
the homogeneous coordinate of  $P_0$
is  $(1,0,...,0)$  
after coordinate transformation. Let $x=(x_0, x_1, ..., x_n)$
be the homogeneous  coordinates of ${\mathbb{P}}^n$.

 The idea of the proof is the following.
Let $V$  be the vector space of homogeneous polynomials 
of degree $a$  passing  through the point $P_0$. For every closed
point $x$, we will  
construct a map $\xi_x $ from $V$  to 
${\mathcal{O}}_{x, X}/{\mathcal{M}}_{x,X}^2$
such that $\xi_x$  is surjective 
for all closed points $x\neq P_0$  and 
$\xi_{P_0}$  is 
surjective from $V$  to 
${\mathcal{M}}_{P_0,X}/{\mathcal{M}}_{P_0,X}^2$.
Let $S_x$  be the set of smooth hypersurfaces 
$H$  in  $V$ such that $x$  is a singular point of  
$H\cap  X$ or  $X\subset  H$.  Let $S$  be the set of closed points of
a closed subset of  projective variety $X\times V$:

$$S=\{<x, H>| x\in  X, H\in  S_x\}.$$

 Let
$p_2: S\rightarrow V$ be the projection. We will show that
the image $p_2(S)$ is a closed subset of 
$V$. So a general member of $V$  intersects $X$
with a smooth  subvariety of codimension 1.
By standard vanishing theorems, we will obtain the  
irreducibility.

For the simplicity, we will first give detail when 
the degree is 2. Higher degree case can be proved in Step 6
by the same method.\\

{\bf{Step 1.}} Let  $V$  be the vector space of the hypersurfaces
of degree 2 passing through $P_0$, then 
a general member of $V$ is smooth.

Let $H$  be a    hypersurface  defined by a homogeneous 
polynomial  $h$   of   degree   2
passing  through  $P_0$, then 
$$h=\sum_{j=1}^na_{0j}x_0x_j
+\sum_{i=1}^n\sum_{j=i}^n a_{ij}x_ix_j.
$$  
Since $\frac{\partial h}{\partial x_j}(P_0)=a_{0j}$, 
$H$  is nonsingular  at 
$P_0$  if at least one $a_{0j}\neq  0$.

The dimension of  $V$  as a vector space  is  
$$\mbox{dim}_k V=
\frac{(n+2)(n+1)}{2}-1=\frac{n^2+3n}{2}.$$

By  Euler's formula, the hypersurface $H$ is singular at a point
 $x=(x_0, x_1,..., x_n)$ if and only if 
$$\frac{\partial h}{\partial x_0} =\frac{\partial h}{\partial x_1} =...=
\frac{\partial h}{\partial x_n}=0. 
$$
It is    a  system of  linear  equations

     $$\quad\quad\quad \quad  a_{01}x_1 + a_{02}x_2 + \cdots    + a_{0n}x_n  =0$$
    $$ a_{01}x_0 + 2a_{11}x_1 + a_{12}x_2 + \cdots    +  a_{1n}x_n  =0$$
    $$\cdot\cdot\cdot\cdot\cdot\cdot\cdot\cdot\cdot\cdot\cdot\cdot\cdot
\cdot\cdot\cdot\cdot\cdot\cdot\cdot\cdot\cdot
\cdot\cdot\cdot\cdot\cdot\cdot\cdot\cdot\cdot\cdot\cdot\cdot$$ 
    $$a_{0n}x_0       + a_{1n}x_1 +a_{2n}x_2 +\cdots    +  2a_{nn}x_n  =0$$

The above system has a solution in  ${\mathbb{P}}^n$
if and only if  the determinant   of  
the following 
 symmetric  matrix  $A$ is zero,

\[\left( \begin{array}{llll}
 0  & a_{01}  & a_{02} \cdots & a_{0n}\\
a_{01} & 2a_{11} & a_{12} \cdots & a_{1n}\\
\multicolumn{4}{c}\dotfill\\
a_{0n} & a_{1n}  & a_{2n} \cdots  & 2a_{nn}
\end{array}   \right). \]

Considering  $(a_{01}, a_{02},..., a_{(n-1)n}$
as a point in the projective space  ${\mathbb{P}}^{\frac{n^2+3n}{2}-1}$,
the system has solutions  only on the hypersurface  defined  by
det$A=0$.  So  the degree 2  hypersurface $H$ in $V$  is nonsingular  on an 
open subset of  ${\mathbb{P}}^{\frac{n^2+3n}{2}-1}$, i.e., a general 
member  $H$   of   $V$  is smooth. Thus among the hypersurfaces of degree 
2 passing through $P_0$, a general member is smooth.\\

{\bf{Step 2.}} There is a map $\xi_x$  from
$V$  to ${\mathcal{O}}_{x, X}/{\mathcal{M}}_{x,X}^2$
such that $\xi_x$  is surjective 
for all closed points $x\neq P_0$  and 
$\xi_{P_0}$  is 
surjective from $V$  to 
${\mathcal{M}}_{P_0,X}/{\mathcal{M}}_{P_0,X}^2$.

Let  $x$  be a   closed  point of  $X$  and  define  $S_x$  to  be the set of
smooth hypersurfaces  $H$  (defined by  $h$) of degree
2 in $V$
such that      $x$  is  a  singular 
point of  $H\cap  X$ ($\neq X$) or $X\subset H$.   Fix  a 
 hypersurface  $H_0$
of degree 2 in $V$ such  that  $x$  is not a point of  $H_0$. Let $h_0$
be the  defining   homogeneous polynomial of $H_0$, then
$h/h_0$  gives  a  regular function 
 on  ${\mathbb{P}}^n-H_0$.  
When restricted     to   $X$, it is  a regular 
function on  $X-X\cap  H_0$.

Let ${\mathcal{M}}_{x,X}$ be the maximal ideal of the 
local ring ${\mathcal{O}}_{x, X}$ at $x$. 
  Define   a map $\xi_x$   from  the  vector  space   $V$   to 
${\mathcal{O}}_{x, X}/{\mathcal{M}}_{x,X}^2$    as follows:
for every  element $h$  in  $V$   (a homogeneous  polynomial 
of degree 2 such that the  corresponding  hypersurface  $H$ is smooth and 
passes through the fixed point $P_0$),  the image  $\xi_x(h)$   is the image  
of  $h/h_0$   in the local  ring  
${\mathcal{O}}_{x, X}$  
 modulo   ${\mathcal{M}}_{x,X}^2$.   It is easy to see that
$x$  is a point of   $H\cap   X$  if and only if  the 
image $\xi_x(h)$ of the defining polynomial  $h$ of  $H$ 
is  contained in   ${\mathcal{M}}_{x,X}$.  And
$x$ is singular   on  $H\cap   X$  if and only if  the 
image $\xi_x(h)$   is  contained in   ${\mathcal{M}}_{x,X}^2$, 
because  the local ring  ${\mathcal{O}}_{x,X}/\xi_x(h)$
will not be regular. 
So there is the following one-to-one  correspondence
$$H\in  S_x \Longleftrightarrow  h\in   ker  \xi_x.  
$$      
Since   $x$ is a closed point and  the  ground field is  
algebraically closed of characteristic 0,  the maximal  ideal  
${\mathcal{M}}_{x,X}$    is   generated   by  linear  forms  in  the coordinates.
Let  $d$  be the dimension of   $X$,  then  the vector space  
${\mathcal{O}}_{x, X}/{\mathcal{M}}_{x,X}^2$  has dimension  $d+1$  
over  $k$.

We will show that
the map  $\xi_x$  is surjective from
$V$  to 
${\mathcal{O}}_{x, X}/{\mathcal{M}}_{x,X}^2$
 if $x\neq  P_0$  and $\xi_{P_0}$ is surjective 
 from $V$  to ${\mathcal{M}}_{P_0,X}/{\mathcal{M}}_{P_0,X}^2$.

Let $U_i=\{(x_0,..., x_n)\in {\mathbb{P}}^n| x_i\neq 0\}$.
Then $\{U_0,..., U_n\}$  is an affine open cover of 
${\mathbb{P}}^n$. In $U_1$, 
we choose the local coordinates in the following 
 $$y_1=\frac{x_0}{x_1}, \quad y_2=\frac{x_2}{x_1}, \quad...,  \quad y_n=\frac{x_n}{x_1}.$$
Let $P$ be a closed  point in $X\cap U_1$
and $(a_1,..., a_n)$  be  the local  coordinate
of  $P$  in $U_1$. 
 We choose $h_0$ to be
$x_1^2$, then
$$
 \frac{h}{h_0}= a_{01}(\frac{x_0}{x_1})+ a_{11}+a_{12}(\frac{x_2}{x_1})+...+
a_{1n}(\frac{x_n}{x_1})+ other\quad terms.
$$
$$
\quad\quad\quad\quad\quad=a_{01}(y_1-a_1)+a_{12}(y_2-a_2)+...+a_{1n}(y_n-a_n)
+c+other\quad terms,
 $$
 
 where the constant 
 $$c=a_{11}+a_{01}a_1+a_{12}a_2+...+a_{1n}a_n.$$ 
 From the expression of $ \frac{h}{h_0}$, we 
 know that 
 $$\{h/h_0|h\in  V\}\rightarrow {\mathcal{O}}_{P, {\mathbb{P}}^n}/{\mathcal{M}}_{P, {\mathbb{P}}^n}^2
 $$
 is surjective.
 So
 $\xi_P$  is surjective  (\cite{H}, page 32).
 
 For any closed  point $P=(a_1,..., a_n)$
 in $U_i\cap X$, $i=2,..., n$, choose $y_j=x_j/x_i$ as local coordinates, 
 similar calculation shows that $\xi_P$  is surjective.

In $U_0$, let $h_0=x_0^2$  and 
$y_i=x_i/x_0$,   then
$$\frac{h}{h_0}=\sum_{j=1}^na_{0j}y_j+\sum_{i=1}^n\sum_{j=i}^n a_{ij}y_iy_j.
$$
There is no constant term in $h/h_0$,
 so
the map  $\xi_{P_0}$  is not surjective to 
${\mathcal{O}}_{P_0, X}/{\mathcal{M}}_{P_0, X}^2$ 
but    
surjective  to  
${\mathcal{M}}_{P_0, X}/{\mathcal{M}}^2_{P_0, X}$.

If $P=(a_1,..., a_n)\neq P_0=(0,...0)$ is a closed  point
in $U_0\cap X$, then write  $y_i=(y_i-a_i)+a_i$, we have 

$$\frac{h}{h_0}= \sum_{j=1}^na_{0j}(y_j-a_j)+\sum_{j=1}^na_{0j}a_j+\sum_{i=1}^n\sum_{j=i}^n a_{ij}[(y_i-a_i)+a_i][(y_j-a_j)+a_j]
$$
$$
=c+I+II,
$$

where the constant term

$$c=\sum_{j=1}^na_{0j}a_j+\sum_{i=1}^n\sum_{j=i}^n a_{ij}a_ia_j,
$$

the linear term with respect to  $y_i-a_i$ is complete

$$ I=\sum_{j=1}^na_{0j}(y_j-a_j)+\sum_{i=1}^n\sum_{j=i}^n a_{ij}[a_i(y_j-a_j)+a_j(y_i-a_i)],
$$

 and the degree 2 term with respect to $y_i-a_i$
 
 $$ II=\sum_{i=1}^n\sum_{j=i}^n a_{ij}(y_i-a_i)(y_j-a_j).
 $$

 Since  $P=(a_1,..., a_n)\neq  (0,..., 0)$,
 there is at least one  $i$, such that $a_i\neq 0$, 
 $1\leq i\leq n$.  So the arbitrary constant 
 $a_{ii}a_i^2$  is a term in 
 $c$. 
 The above expressions of  constant $c$  and linear term 
 $I$ show that on $U_0$, $\xi_x$  is surjective if 
 the closed  point    $x\neq P_0$.

Considering  the kernel of  the  map  
$$\xi_x: V\longrightarrow    {\mathcal{O}}_{x, X}/{\mathcal{M}}_{x,X}^2,
$$
if  $x\neq  P_0$, 
the kernel  as a vector space has dimension 

$$  \mbox{dim}_{\mathbb{C}}ker \xi_x=\frac{n(n+3)}{2}-d-1. $$

Therefore   $S_x$  is  a  linear  system  of
hypersurfaces    with
dimension   $\frac{n(n+3)}{2}-d-2$  if 
$x\neq  P_0$.  If  $x=P_0$, then  the  projective 
dimension
of  $S_{P_0}$  is   
 $\frac{n(n+3)}{2}-d-1$.\\

{\bf{Step 3.}}
If  $V$ is considered   as a  projective space, 
then  $X\times  V$ is a projective  variety.  Let the subset  
$S\subset  X\times  V $   consist    of all   pairs   $<x, H>$
such that  $x\in  X$  is  a  closed  point  and  
$H\in   S_x$. Then the dimension of 
$S$  is less than the dimension of 
$V$.

        $S$  is the set of  closed  points   of  a closed  subset   of
$  X\times  V $  and  we   give    a reduced  induced  scheme structure
to  $S$.  The  first  projection   
$p_1:  S\rightarrow   X$  is surjective. 
If $x\neq  P_0$,  
 the  fiber 
 $p_1^{-1}(x)$   
is a projective space   with  dimension   
$\frac{n(n+3)}{2}-d-2.$ 
The special fiber  $p_1^{-1}(P_0)$
is a projective space with
 dimension $\frac{n(n+3)}{2}-d-1.$

Let  $S=\cup_{i=0}^mS_i $  be
 an  irreducible  decomposition.
 Then  every  $p_1(S_i)$  
 is closed  and 
 there is an $i$, such that 
 $p_1(S_i)=X$.  For every  $S_i$   with
   $p_1(S_i)=X$,  there is an open subset 
   $U_i\subset S_i$ such that for every $x\in U_i$,
   the fiber $p_1^{-1}(x)$ has  constant dimension 
   $n_i$
   (\cite{S}, Chapter 1, Section 6.3, Theorem 7).  Let $x\in \cap U_i$, since  
   the fiber  $p_1^{-1}(x)$  is irreducible, it is contained in some $S_i$. Suppose $p_1^{-1}(x)\in S_1$.
   Let $f_1$ be the restriction of  $p_1$ on 
   $S_1$, i.e., $p_1|_{S_1}=f_1$, then  
   $p_1^{-1}(x)\subset  f_1^{-1}(x)$  since  $p_1^{-1}(x)$
   is  irreducible. 
   The opposite inclusion  is  obvious, so
   $p_1^{-1}(x)=f_1^{-1}(x)$ for $x\in \cap U_i$
   and  $n_1=\frac{n(n+3)}{2}-d-2$. 
   
   Since  $f_1$  is surjective and $S_1$  is one
   irreducible  component  of  $S$, for every 
   $x\in X$, the fiber  $f_1^{-1}(x)$  
   is  not empty and contained in 
    $p_1^{-1}(x)$.  But  the dimension of 
   $f_1^{-1}(x)$  is  at least 
   $\frac{n(n+3)}{2}-d-2$, so 
   for every 
   $x\in X$, $x\neq P_0$,  $p_1^{-1}(x)=f_1^{-1}(x)$. 
   
   Hence  $S_1$     has    
dimension   (\cite{S}, Chapter 1, Section 6.3)

$$ [\frac{n(n+3)}{2}-d-2]+d=\frac{n(n+3)}{2}-2.$$

 If there is a component  $S_j$
 in $S$  such that $p_1(S_j)\neq X$, then
 the dimension of $S_j$ is not
 greater than 
 the dimension of $S_1$
 (\cite{S}, Chapter I, Section 6.3, Theorem 7). So if 
 $S$  is not irreducible,
 then for all components $S_i$ in $S$,
 $S_1$  has the maximum dimension
 $ [\frac{n(n+3)}{2}-d-2]+d=\frac{n(n+3)}{2}-2$, which is the 
 dimension of $S$. \\

 {\bf{Step 4.} }
 A general member $H$ of $V$
 intersects $X$  with a smooth codimension 1
 subvariety on $X$. 

Looking   at  the   second  projection    (a proper  morphism )
$p_2: S\rightarrow    V$.  The dimension of the image

$$  \mbox{dim}  p_2(S)\leq  \mbox{dim}  S=\frac{n(n+3)}{2}-2.$$

Since  $S$  is closed   in   $  X\times  V $
and the dimension of $V$ (as a projective space)
is  $\frac{n(n+3)}{2}-1$,
$V-p_2(S)$  is an open   subset   of  
$V$. This implies that a general member $H$ of  $V$
intersects  $X$  
with a smooth  variety $X\cap  H$.\\

{\bf{Step 5.} }
  $X\cap  H$
is   irreducible.

From the short exact  sequence
$$ 0\longrightarrow 
 {\mathcal{O}}_{{\mathbb{P}}^n}(-H)
\longrightarrow 
 {\mathcal{O}}_{{\mathbb{P}}^n}
\longrightarrow 
 {\mathcal{O}}_{H}
\longrightarrow 
0,
$$
since  
$H^1({\mathbb{P}}^n, {\mathcal{O}}_{{\mathbb{P}}^n}(-H))=0$  
(\cite{H}, Page 225, Theorem 5.1),
we have   a  surjective map 
$$H^0({\mathbb{P}}^n,{\mathcal{O}}_{{\mathbb{P}}^n})=k
\longrightarrow  H^0(H, {\mathcal{O}}_{H}).$$
 $H^0(H, {\mathcal{O}}_{H})=k$   implies that 
 the hypersurface  $H$  is connected.

Since  $X$  is  closed in  ${{\mathbb{P}}^n}$, 
$H|_X$  is ample  on  $X$.   By 
Kodaira  Vanishing  Theorem,
$H^1(X, {\mathcal{O}}_{X}(-H))=0$  (\cite{KM}, Page  62).
Applying the short exact sequence
$$ 0\longrightarrow 
 {\mathcal{O}}_X(-H)
\longrightarrow 
 {\mathcal{O}}_X
\longrightarrow 
 {\mathcal{O}}_{H\cap  X}
\longrightarrow 
0,
$$ 
we  get  
$$H^0(H\cap X,  {\mathcal{O}}_{H\cap  X})=H^0(X, {\mathcal{O}}_X)=k.$$
Thus  the intersection  $H\cap  X$  is connected.  
Therefore  for a general hypersurface 
$H$  of degree 2,   $H\cap  X$  is   smooth  and  irreducible.

We have proved that  a general  smooth  hypersurface   of  degree  2  
 passing   through   
$P_0$   intersects   $X$  with  an  irreducible  
 smooth  subvariety   of codimension  1.

 {\bf{Step 6.} } Degree $a>2$ case.

 Let $W$  be the vector space of 
 hypersurfaces $H$ of degree $a>2$
 such that $P_0\in  H$. Then any element of
 $W$  can be written in the following form
 $$g=c_0x_0^{a-2}h+c_1x_1^{a-2}h+...+c_nx_n^{a-2}h
 +other \quad terms,
 $$
 where 
 $$h=\sum_{j=1}^na_{0j}x_0x_j
+\sum_{i=1}^n\sum_{j=i}^n a_{ij}x_ix_j,
$$  
is the hypersurface   of  degree  2 in Step 1. 

It is easy to see from the calculation of Step 2
that in  each affine open subset $U_i$, 
we have
$$\{\frac{h}{x_i^2}| h\in V\}\subset  
\{\frac{g}{x_i^a}| g\in W\}. 
$$ 
 So again the map $\xi_x$ from $W$
 to  ${\mathcal{O}}_{x, X}/{\mathcal{M}}_{x,X}^2$
 is surjective and 
 $\xi_{P_0}$ is surjective from 
 $W$  to ${\mathcal{M}}_{P_0,X}/{\mathcal{M}}_{P_0,X}^2$.

 For any degree $a>2$  hypersurface, by counting the dimension correctly as above,
 we can similarly show that  a general 
 hypersurface passing through $P_0$ 
 intersects $X$  with an irreducible smooth projective
 variety of codimension 1 on $X$. In fact, the vector space 
 $W$
  has dimension greater than $(n^2+3n)/2$,
  and dimension of $p_2(S)$  is less than the dimension of 
  $W$ as a projective space.
   So the whole argument works.

\begin{flushright}
 Q.E.D. 
\end{flushright}

\begin{remark}
We only need the intersection part $X\cap H$  is irreducible and smooth.
Outside $X$, $H$  being  smooth or not  does not play any role.  
\end{remark}

\begin{remark}  From the proof, we see that if 
the point $P_0$  is a point outside  $X$, the theorem still
holds since the map $\xi_x$  is surjective for all 
$x\in X$. 

\end{remark}

\begin{remark}
Let $L$ be the vector  space of hyperplanes passing through 
$P_0$, then its dimension is $n$. If dim$X=d=n-1$, then the map $\xi _x$ may not be surjective. So the above  proof
does not work  for hyperplanes.
For  general hyperplane sections, Theorem 2.1 is still true
but the proof is different.
\end{remark}

The proof of the following Theorem 2.5 has been communicated 
 to me by 
Igor Dolgachev.

\begin{theorem}
Let $X$  be an irreducible smooth projective variety
of dimension at least 2 in ${\mathbb{P}}^n$. 
Let $P_0$ be a closed point on 
$X$. Then a general hyperplane  passing through 
$P_0$  is irreducible  and smooth. 

\end{theorem}

$Proof.$ If $X$ is a hyperplane, then the theorem is true.
So we may assume that 
$X$  is not a hyperplane.

Consider the dual projective space ${\mathbb{P}}^{n*}$. The 
hyperplanes passing through $P_0$  in 
${\mathbb{P}}^n$ is a hyperplane 
  $H^*$  in ${\mathbb{P}}^{n*}$.  
  A hyperplane  $B$ intersects 
  $X$  with  a  singular  point  $x\in B\cap  X$
  if  $B$  is tangent to  $X$  at $x$. 
  Let 
  $X^*$  be the dual space of $X$, i.e.,
  the set of hyperplanes tangent to $X$
  at some point. Then the dimension 
  of $X^*\leq n-1$  (\cite{AG}, Section 2.5, 3.1). 
  If $X^*=H^*$, then $(X^*)^*=X=H$, which is not 
  possible by our assumption. So
  $X^*\cap H^*$  has dimension at most
  $n-2$   and $H^*-X^*$ is an open subset of
  $H^*$. Any point of $H^*$  away from 
   $X^*$  corresponds   to a hyperplane 
   in ${\mathbb{P}}^n$   which cuts $X$  with 
   a smooth subvariety of codimension 1 on $X$.
   By   Kodaira  Vanishing Theorem, 
   the intersection subvariety is irreducible.

\begin{flushright}
 Q.E.D. 
\end{flushright}

\begin{remark}
By the proof, we see that if  the point $P_0$    is  not
 a point on $X$,  then Theorem 2.5  is  still  true.
\end{remark}

\begin{remark}

By Veronese embedding and considering the dual variety, we can prove Theorem 2.1 in a
much easier  and geometric way. The advantage 
of the long proof is that  the idea  can be 
used  to the case of $q$ points, $q>1$. 
\end{remark}

\section{Hypersurface Sections Passing Through $q$ Points}

\begin{theorem}
If  $X$  is an irreducible   smooth projective  variety
of dimension at least  2 in ${\mathbb{P}}^n$, then  for any $q\leq n+1$
closed points 
 $P_0,$ $P_1,..., P_{q-1}$  on  $X$
 in a general position
and any degree $a>1$,
a general    hypersurface 
$H$  of  degree  $a$ passing through  these $q$  points 
intersects $X$  with 
  an irreducible  smooth  codimension 1  subvariety  
on  $X$.
\end{theorem}

$Proof$.   By  the proof of  Theorem 2.1, Step 5,
for any hypersurface 
$H$, $H\cap X$  is connected. So we only need to show that 
$H\cap X$  is nonsingular for a general  $H$  passing through
these $q$  points.  

{\bf{Step 1.} } There is a hyperplane $H$ in 
${\mathbb{P}}^n$ such that $H$ does not contain any   point   $P_s$,
$s=0,1,..., (q-1)$.

 Let $(x_0, x_1,..., x_n)$ be the homogeneous coordinates
 of ${\mathbb{P}}^n$, and $a_0x_0+a_1x_1+...+a_nx_n=0$
 be the hyperplane $H$.
 
 Let $P_s=(b_{s0}, b_{s1},..., b_{sn})$. Consider
 $(a_0, a_1,..., a_n)$ such that 
 $a_0b_{s0}+a_1b_{s1}+...+ a_n b_{sn}=0$. We may look at
 it as a hyperplane $H^*_s$ in the dual projective space 
 ${\mathbb{P}}^{n*}$.   Choose $c_0, ..., c_n$ in the open 
 subset ${\mathbb{P}}^{n*}\setminus  \cup _{s=0}^{q-1} H_s^*$,  then 
 for every $P_s$, $c_0b_{s0}+...+c_nb_{sn}\neq 0$. 
So there is a hyperplane $H$ defined by 
$c_0x_0+...+c_nx_n=0$  that  does not contain any  point
$P_s$, $s=0,1,..., q-1$.  

{\bf{Step 2.} } We may change the coordinates such that 
$P_0=(1, 0, ..., 0)$ and $H$ is defined by 
$x_0=0$ such that $H$ does not contain any 
$P_s$, $s=0,1,..., q-1$.

By Step 1, we may choose hyperplane $H$: $c_0x_0+...+c_nx_n=0$ 
such  that  $H$   does not  contain  any  point  $P_s$,
$s=0, 1,..., q-1$. Define 
the new coordinates 
$$X_0=c_0x_0+...+c_nx_n,$$
 $$X_1=\sum_{j=0}^na_{1j}x_j, $$
 $$......$$
$$X_n=\sum_{j=0}^na_{nj}x_j,$$
 where the coefficients $a_{ij}$  satisfy the following
 system of linear equations
 $$\sum_{j=0}^na_{1j}b_{1j}=0,
 $$
$$\sum_{j=0}^na_{2j}b_{2j}=0,
$$
$$......,$$
$$\sum_{j=0}^na_{nj}b_{nj}=0.
$$
The $n$  points $P_1,..., P_n$   being linearly independent
guarantees  that the above  linear  transformation   is  well-defined. 
Since  $c_0b_{00}+c_1b_{01}+...+ c_nb_{0n}\neq  0$
 the new coordinate of $P_0$ is $(1, 0, ..., 0)$
and the plane $H: X_0=0$ does not contain 
 any 
$P_s$, $s=0,1,..., q-1$.

{\bf{Step 3.} } Let  $V$  be the vector space of the hypersurfaces
of degree 2 passing through these $q\leq n+1$ closed points $P_0,...P_{q-1}$
in general position, then the map
$\xi_x$  from
$V$  to ${\mathcal{O}}_{x, X}/{\mathcal{M}}_{x,X}^2$
defined  in  the proof  of  Theorem 2.1, Step 2, 
   is surjective 
for all closed points $x\in U_0$, $x\neq P_s$   and 
$\xi_{P_s}$  is 
surjective from $V$  to 
${\mathcal{M}}_{P_s,X}/{\mathcal{M}}_{P_s,X}^2$,
$s=0,1,..., q-1$, where 
${\mathcal{M}}_{P_s,X}$  is  the  maximal  ideal  of
${\mathcal{O}}_{x, X}$.

By Step 2, we may assume that $P_0=(1,0,...,0)$
and the hyperplane $H_0$ defined by $x_0=0$ does not 
contain any point $P_s$, $s=0,1,..., q-1$, $q\leq n+1 $.
By Step 1, proof of Theorem 2.1,
a homogeneous 
polynomial  $h$   of   degree   2
passing  through  $P_0$ is of the form
$$h=\sum_{j=1}^na_{0j}x_0x_j
+\sum_{i=1}^n\sum_{j=i}^n a_{ij}x_ix_j.
$$   
In $U_0={\mathbb{P}}^n-H_0$, let $h_0=x_0^2$  and 
$y_i=x_i/x_0$,   then
$$\frac{h}{h_0}=\sum_{j=1}^na_{0j}y_j+\sum_{i=1}^n\sum_{j=i}^n a_{ij}y_iy_j.
$$
Let $P_s=(b_{s1}, b_{s2}, ..., b_{sn} )$ in $U_0\cong {\mathbb{A}}^n$.
Since $h(P_s)=0$ but $h_0(P_s)\neq 0$, we have the following 
$q-1$ equations
$$\frac{h}{h_0}(P_s)=\sum_{j=1}^na_{0j}b_{sj}+\sum_{i=1}^n\sum_{j=i}^n a_{ij}b_{si}b_{sj}=0,
$$
for $s=1,2, ..., q-1.$ This is a system of $(q-1)\leq n$ linear equations 
with $(n^2+3n)/2$  variables   $a_{ij}$.

The coefficient matrix $A$ of  this linear system is 

\[\left( \begin{array}{llllll}
 b_{11}  & b_{12}\quad \cdots  & b_{1n} & b_{11}^2 & b_{11}b_{12}\quad\quad \cdots & b_{1n}^2\\
 b_{21}  & b_{22}\quad\cdots  & b_{2n} & b_{21}^2 & b_{21}b_{22}\quad\quad \cdots & b_{2n}^2\\
\multicolumn{6}{c}\dotfill\\
 b_{(q-1)1}  & b_{(q-1)2}\cdots  & b_{(q-1)n} & b_{(q-1)1}^2 & b_{(q-1)1}b_{(q-1)2}\cdots & b_{(q-1)n}^2\\
\end{array}   \right). \]
\\

The matrix $B$ of the  first $n$ columns
and all  $(q-1)$ rows is  

\[\left( \begin{array}{lll}
 b_{11}  & b_{12}\quad \cdots  & b_{1n} \\
 b_{21}  & b_{22}\quad\cdots  & b_{2n} \\
\multicolumn{3}{c}\dotfill\\
 b_{(q-1)1}  & b_{(q-1)2}\cdots  & b_{(q-1)n} \\
\end{array}   \right). \]
\\

Because the $q-1$ points $P_1, ..., P_{q-1}$ are 
in general position, the rank of 
the above matrix 
$B$  is $q-1\leq n$. 
So the  system defined by $A$ is consistent
and each $a_{0j}$ varies independently. The dimension of 
solutions of the system $\frac{h}{h_0}(P_s)=0$,
$s=1,..., q-1$,  
is 
$$\frac{n^2+3n}{2}-(q-1)$$ 
as a vector space. 

For any closed point $P=(a_1, a_2,..., a_n)$
in $U_0$, we have 

$$\frac{h}{h_0}= \sum_{j=1}^na_{0j}(y_j-a_j)+\sum_{j=1}^na_{0j}a_j+\sum_{i=1}^n\sum_{j=i}^n a_{ij}[(y_i-a_i)+a_i][(y_j-a_j)+a_j]
$$
$$
=c+I+II,
$$

where the constant term

$$c=\sum_{j=1}^na_{0j}a_j+\sum_{i=1}^n\sum_{j=i}^n a_{ij}a_ia_j,$$

the linear term with respect to  $y_i-a_i$

$$ I=\sum_{j=1}^na_{0j}(y_j-a_j)+\sum_{i=1}^n\sum_{j=i}^n a_{ij}[a_i(y_j-a_j)+a_j(y_i-a_i)],
$$

 and the degree 2 term with respect to $y_i-a_i$
 
 $$ II=\sum_{i=1}^n\sum_{j=i}^n a_{ij}(y_i-a_i)(y_j-a_j).
 $$

If $P\neq P_0$, then there is at least one $i$, such that
$a_i\neq 0$,  $i=1,..., n$. So the constant $c$ can be any number since it has 
a term $a_{ii}a_i^2$, where $a_{ii}$ varies independently and 
$a_i\neq 0$. We already see that $a_{01}, a_{02}, ..., a_{0n}$
vary independently,
so  $c+I$  is a  complete linear 
form with every linear term.
Thus we have a surjective map from
$ V$  to  ${\mathcal{O}}_{P, {\mathbb{P}}^n}/{\mathcal{M}}_{P,{\mathbb{P}}^n }^2$. 
Since there is a natural surjective map from 
${\mathcal{O}}_{P, {\mathbb{P}}^n}/{\mathcal{M}}_{P, {\mathbb{P}}^n}^2$
to ${\mathcal{O}}_{P, X}/{\mathcal{M}}_{P, X}^2$
(\cite{H}, Page 32), 
by the above expression of $h/h_0$, 
if the closed point $P=P_s$, $s=1,2,..., q-1$, then $\xi_{P_s}$
is surjective from 
$ V$  to ${\mathcal{O}}_{P_s, {\mathbb{P}}^n}/{\mathcal{M}}_{P_s, {\mathbb{P}}^n}^2$
so surjective to ${\mathcal{O}}_{P_s, X}/{\mathcal{M}}_{P_s, X}^2$.

{\bf{Step 4.} } The map
$\xi_x$  from
$V$  to ${\mathcal{O}}_{x, X}/{\mathcal{M}}_{x,X}^2$
   is surjective 
for all closed points $x\in H_0$.

Since on $H_0$, $x_0=0$, every element in
$V$  is of the same form
 $$h=\sum_{j=1}^na_{0j}x_0x_j+ \sum_{i=1}^n\sum_{j=i}^n a_{ij}x_ix_j.
 $$
Let $P=(0, a_1, a_2,..., a_n)$ be a closed point
on $X\cap H_0$, then 
at least one $a_i\neq 0$. Suppose that $a_1\neq 0$. Let $h_1=x_1^2$,
define $y_i=x_i/x_1$, 
 then 
$$\frac{h}{h_1}= a_{01}y_0+ a_{02}y_0y_2+...+a_{0n}y_0y_n+ a_{11}+a_{12}y_2+...+a_{1n}y_n+\sum_{i=2}^n\sum_{j=i}^n a_{ij}
y_iy_j.
$$
We can rewrite it into three parts
$$
\frac{h}{h_1} =c+I+II,
$$
where the constant 
$$c= a_{11}+ \sum_{i=2}^n	a_{1i}a_i+\sum_{i=2}^n\sum_{j=i}^n a_{ij}a_ia_j,
$$
the linear term
$$I=a_{01}y_0+ \sum_{i=2}^n	a_{1i}(y_i-a_i)+
\sum_{i=2}^n\sum_{j=i}^n a_{ij}[a_j(y_i-a_i)+a_i(y_j-a_j)],
$$
and the degree 2 term
$$II= \sum_{i=2}^na_{0i}y_0(y_i-a_i)+ \sum_{i=2}^n\sum_{j=i}^n a_{ij}(y_i-a_i)(y_j-a_j).
$$
Since $h/h_1$   contains 
the  complete 
linear  form  with every linear term, 
 the map from $V$ to
${\mathcal{O}}_{P, X}/{\mathcal{M}}_{P,X}^2$
 at any point 
$P=(0, a_1,..., a_n)$ on $H_0$ with $a_1\neq 0$ is 
surjective. In general,  if 
$a_i\neq 0$, choose $h_i=x_i^2$, then 
by the same calculation, we can show that 
the map $\xi_P$ is surjective from 
$V$  to  
${\mathcal{O}}_{P, X}/{\mathcal{M}}_{P,X}^2$.

{\bf{Step 5.} } Let  $x$  be a   closed  point of  $X$  and  define  $S_x$  to  be the set of
smooth hypersurfaces  $H$  (defined by  $h$) of degree
2 in $V$
such that      $x$  is  a  singular 
point of  $H\cap  X$ ($\neq X$) or $X\subset H$. 
If  $V$ is considered   as a  projective space, let the subset  
$S\subset  X\times  V $   consist    of all   pairs   $<x, H>$
such that  $x\in  X$  is  a  closed  point  and  
$H\in   S_x$. Then the dimension of 
$S$  is less than the dimension of 
$V$.

As a projective space, $V$  has dimension 
$$\frac{n^2+3n}{2}- (q-1)-1=\frac{n^2+3n}{2}-q
\geq 
 \frac{n^2+n}{2}-1,$$ 
where $q\leq n+1$. Let  $d$
be the dimension of of  $X$. 
 Consider the first 
projection $p_1: S\rightarrow X$,  for all closed points
 $x\in X$, $x\neq P_s$, $s=0,1,...q-1$,
 as a  vector space, ${\mathcal{O}}_{x, X}/{\mathcal{M}}_{x,X}^2$
 has dimension $d+1$. 
  By Step 3 and 4,
 the dimension of the fiber $p_1^{-1}(x)$ is 
 the dimension of  
the  kernel  of the map  
$\xi_x$, so  as a projective space, 
 $${\mbox{dim}}p_1^{-1}(x)= 
 (\frac{n^2+3n}{2}-q)-(d+1)= 
 \frac{n^2+3n}{2}-q-d-1.$$

 The vector space ${\mathcal{M}}_{P_s, X}/{\mathcal{M}}_{P_s,X}^2$
 has dimension  $d$.  
 The dimension of the fiber over the  point  $P_s$ is
 $${\mbox{dim}}p_1^{-1}(P_s)=(\frac{n^2+3n}{2}-q)-d= \frac{n^2+3n}{2}-q-d.
 $$
 
The dimension of $S$  is  
$$\frac{n^2+3n}{2}-q-d-1+ d=\frac{n^2+n}{2}-q-1.
$$ 
So 
$${\mbox{dim}}(V)-{\mbox{dim}}(S)=(n^2+3n)/2- q)-(\frac{n^2+3n}{2}-q-1)=1.
$$
Let $p_2: S\rightarrow V$ be the second projection,
then 
$${\mbox{dim}}(p_2(S))\leq {\mbox{dim}}(S)< {\mbox{dim}}(V).$$

 Since  $p_2(S)$  is a closed subset of $V$, a general member 
of $V$  intersects $X$ with a smooth subvariety of codimension 1 on $X$. 
By Step 5 of proof of Theorem 2.1, 
it is also connected so irreducible. This
proves the degree 2 case.

{\bf{Step 6.} } Let $h$ be an element of $V$ in Step 3. Then 
any degree $a>2$ homogeneous polynomial passing through 
the $q$ points $P_0, P_1,..., P_{q-1}$ can be written as
$$f=c_0x_0^{a-2}h+c_1x_1^{a-2}h+...+c_nx_n^{a-2}h+\mbox{other terms},
$$
where  $c_0,..., c_n$  are constants. 
Let $W$ be the set of homogeneous polynomials 
of degree $a$ passing through these $q$ points.
Considering $f/x_i^a$, $i=0,1,..., n$, and
using the same calculation, we can show that 
$\xi_x$  from
$W$  to ${\mathcal{O}}_{x, X}/{\mathcal{M}}_{x,X}^2$
   is surjective 
for all closed points $x\neq P_s$ and is surjective
from $W$  to ${\mathcal{M}}_{P_s, X}/{\mathcal{M}}_{P_s,X}^2$.
Carrying out the dimension calculation as in Step 5, 
we see that the theorem holds for all degree $a>2$.

\begin{flushright}
 Q.E.D. 
\end{flushright}

\section{Linear System of Ample Divisors}

\begin{theorem}
If  $X$  is an irreducible   smooth projective  variety
of dimension   $d\geq 2$, $D$  is an ample  divisor 
on  $X$, then there is an $n_0>0$ such that for all $n\geq n_0$ 
and any closed point $P_0$ on $X$, 
a general member of $|nD|_{P_0}$ is irreducible and  smooth, where $|nD|_{P_0}$
is the linear system of effective divisors in
$|nD|$ passing through 
the point $P_0$.
\end{theorem}

 $Proof.$ Since $D$ is ample, there is an
 $l$ such that the basis
 $\{f_0, ..., f_m\}$ in the vector space 
 $H^0(X, {\mathcal{O}}_X(lD))$
 gives an embedding $\phi$ from  
 $X$ to the projective space 
     ${\mathbb{P}}^m$ by sending a point $x$ on $X$ to
     $(f_0(x), f_1(x), \cdot \cdot\cdot, f_m(x))$ in ${\mathbb{P}}^m$. 
 Let  $W=\phi(X)$, $w_0=\phi(P_0)$. 
 Let $M$  be the set of homogeneous polynomials $h$ of degree
 $a>1$ such that every hypersurface $H$ 
 defined by $h$ contains the point $w_0$  in  
 ${\mathbb{P}}^m$.
 By the proof of the Theorem 2.1, if $w\neq w_0$, then 
 the natural map $\xi_w$ 
 from $M$ to 
     ${\mathcal{O}}_{w,W}/{\mathcal{M}}^2_{w,X}$
     is surjective and 
     $\xi_{w_0}$ is surjective from
     $M$  to ${\mathcal{M}}_{w_0,W}/{\mathcal{M}}^2_{w_0, W}$.

      Let  
     $$ L=\{ f\in H^0({\mathcal{O}}_X(n_0D))| f(P_0)=0\},
     $$
     where $n_0=l^a$, 
     then $\phi^*(M)$  is a subspace of $L$. 
     
     Since
     $\phi$ is an isomorphism 
     between $X$  and $W$,  the  map $\xi_x$ 
    from  $L$  to ${\mathcal{O}}_{x,X}/{\mathcal{M}}_{x,X}^2$
     is also surjective for every point  $x\neq P_0$ in
     $X$ and $\xi_{P_0}$ is surjective from 
     $L$ to ${\mathcal{M}}_{P_0,X}/{\mathcal{M}}^2_{P_0, X}$. 
      The dimension of the kernel of $\xi_x$ as a projective
      space is  dim$L-d-1$
      and the dimension of the kernel of $\xi_{P_0}$
      is  dim$L-d$. Here $L$  is considered as 
      a projective space.

     Because there is a one-to-one correspondence 
     between the  effective 
     divisors in th linear system $|n_0D| $  passing through $P_0$
     and the elements in $L$ up to a nonzero constant, 
     we may think about $L$ as 
     a vector space of effective 
     divisors passing through $P_0$ and linearly equivalent to 
     $n_0D$.

     Let  $x$  be  a  closed  point  of  $X$
     and  
      $S_x$  be the set of effective divisors $E$ in 
$L$  such that $x$  is a singular point of 
$X\cap E$. Let 
$$ S=\{<x, E>|x\in X, E\in S_x \}.
$$

Then $S$  is a closed subvariety of
$X\times L$, where $L$ is the projective space.

A point $x\in X$  is a singular point of 
$f\in L$ if and only if 
the image of $f$ in the natural map 

$$\xi_x: L\rightarrow {\mathcal{O}}_{x,X}/{\mathcal{M}}^2_{x, X}
$$
is zero, i.e., the  image  of  $f$ is in the kernel of 
$\xi_x$.  

Considering the  projections $p_1: S\rightarrow X$
and $p_2: S\rightarrow L$. $p_1$  is surjective.
Every fiber $p_1^{-1}(x)$ is a projective space of 
dimension dim$L-d-1$ if $x\neq P_0$ and $p_1^{-1}(P_0)$
is a projective space of dimension dim$L-d$, where 
 dim$L$ is the dimension of
$L$  as a projective space. By the proof of 
Theorem 2.1, the projective variety $S$   has dimension 
dim$L-d-1+d=$dim$L-1$.

Looking at the second projection $p_2: S\rightarrow L$,
we have dim$p_2(S)\leq $dim$S\leq$ dim$L-1$.
So
$p_2(S)$  is a closed subset of $L$. Thus a  general member of
$L$  is smooth.

Let $E$  be a general smooth element in $L$, then from
$$ 0\longrightarrow 
 {\mathcal{O}}_X(-E)
\longrightarrow 
 {\mathcal{O}}_X
\longrightarrow 
 {\mathcal{O}}_E
\longrightarrow 
0,
$$ 
by Kodaira Vanishing Theorem, we   have 
$$H^0(E,  {\mathcal{O}}_E)
=H^0(X, {\mathcal{O}}_X)=k.$$
Thus  $E$ is connected.

If $n>n_0$, considering  the vector space 
$$L'=\{ f\in H^0({\mathcal{O}}_X(nD))| f(P_0)=0\},
$$
then $L$ is a subspace of $L'$. So 
the map 

$$\xi_x: L'\rightarrow {\mathcal{O}}_{x,X}/{\mathcal{M}}^2_{x, X}$$

is surjective if $x\neq P_0$
and 
$\xi_{P_0}$ is surjective from 
     $L'$ to ${\mathcal{M}}_{P_0,X}/{\mathcal{M}}^2_{P_0, X}$.
The proof still works  for $L'$
by counting the dimensions in the same way. 

\begin{flushright}
 Q.E.D. 
\end{flushright}

\begin{theorem}
If  $X$  is an irreducible   smooth projective  variety
of dimension $d\geq 2$ in ${\mathbb{P}}^n$, $D$  is an ample  divisor 
on  $X$, then there is an $n_0>0$ such that for all $m\geq n_0$ 
and any $q\leq n+1$ closed points $P_0, P_1,..., P_{q-1}$ on $X$
in general position, 
a general member of $|mD|_{q}$ is irreducible and  smooth, 
where $|mD|_{q}$
is the linear system of effective divisors in
$|mD|$ passing through 
these  $q$ points $P_s$, $s=0,1,..., (q-1)$.
\end{theorem}

$Proof$. Let $\phi$ be an embedding as in the proof of Theorem 4.1.
 Let $M$  be the set of homogeneous polynomials $h$ of degree
 $a>1$ such that every hypersurface $H$ 
 defined by $h$ contains the points $\phi(P_s)=w_s$, $s=0,..., q-1$.
 By the proof of the Theorem 3.1, if $w\neq w_s$, then 
 the natural map $\xi_w$ 
 given in the proof  
 from $M$ to 
     ${\mathcal{O}}_{w,W}/{\mathcal{M}}^2_{w,X}$
     is surjective and 
     $\xi_{w_s}$ is surjective from
     $M$  to ${\mathcal{M}}_{w_s,W}/{\mathcal{M}}^2_{w_s, W}$.

      Let  
     $$ L=\{ f\in H^0({\mathcal{O}}_X(n_0D))| f(P_0)=0\},
     $$
     where $n_0=l^a$, 
     then $\phi^*(M)$  is a subspace of $L$. 
     
     Since
     $\phi$ is an isomorphism 
     between $X$  and $W$,  the  map $\xi_x$ 
    from  $L$  to ${\mathcal{O}}_{x,X}/{\mathcal{M}}_{x,X}^2$
     is also surjective for every point  $x\neq P_0$ in
     $X$ and $\xi_{P_0}$ is surjective from 
     $L$ to ${\mathcal{M}}_{P_0,X}/{\mathcal{M}}^2_{P_0, X}$.
     Carrying out the same dimension counting, 
     it is clear that the theorem holds for $|n_0D|_q$. The rest of the proof   is the same as  the proof of  Theorem 4.1.

\begin{flushright}
 Q.E.D. 
\end{flushright}

\section{Linear System of Big Divisors}

Let $X$  be an irreducible normal projective variety
and $D$ a Cartier divisor on $X$.  
If for all $m>0$, $H^0(X, {\mathcal{O}}_X(mD))=0$, then 
the $D$-dimension
$\kappa(D, X)=-\infty$. Otherwise, 
$$\kappa(D, X)=tr.deg_{\mathbb{C}}\oplus_{m\geq 0}
H^0(X, {\mathcal{O}}_X(mD))-1.$$

 If  $h^0(X, {\mathcal{O}}_X(mD))>  0$
for  some  $m\in {\mathbb{Z}}$  and $X$  is normal, 
choose a basis $\{f_0, f_1, \cdot \cdot\cdot, f_n\}$
     of the linear space 
     $H^0(X, {\mathcal{O}}_X(mD))$, it defines a rational 
     map 
     $\Phi _{|mD|}$
     from $X$ to the projective space 
     ${\mathbb{P}}^n$ by sending a point $x$ on $X$ to
     $(f_0(x), f_1(x), \cdot \cdot\cdot, f_n(x))$ in ${\mathbb{P}}^n$. 
We define the  $D$-dimension  (\cite{U}, Definition 5.1),
$$ \kappa (D, X)= \max_m\{\dim (\Phi _{|mD|}(X))\}. 
      $$
      
  From the definition, if $D$  is an effective divisor,
  we have      $0\leq \kappa (D, X)\leq d$,
  where  $d$  is the dimension of  $X$. 
  We say that  
  $D$  is a big divisor if 
$\kappa(D, X)$ is equal to the dimension of 
$X$.

\begin{theorem}Let $X$  be an irreducible smooth
projective variety of dimension $d\geq 2$ and $D$ an effective 
big divisor on 
$X$ such that  $\Phi _{|nD|}$ defines a birational morphism. Then 
there is an open  subset $U$  in $X$  such that 
for any point  
$P_0$  on $U$, a general member of $|nD|_{P_0}$
is smooth, where $|nD|_{P_0}$  is the linear system of effective divisors 
passing through the point $P_0$. 
\end{theorem}

$Proof.$  Since $\Phi_{|nD|}$  is a birational morphism, 
its image $W$ has dimension $d$ and  there is an open subset $U$
in $X$  such that $U$ is isomorphic to $\Phi_{|nD|}(U)$.

 By definition of $D$-dimension, 
  $H^0(X, {\mathcal{O}}_X(nD)$  
  has $d$ algebraically independent nonconstant 
  elements, where  dim$X=d$. 
  Let $(f_0, f_1, ..., f_m)$ be a representation of $\Phi_{|nD|}$.
  Here  $(f_0, f_1, ..., f_m)$ is 
  a basis of 
  $ H^0(X, {\mathcal{O}}_X(nD))$. 
  We may arrange the order such  that  $f_0,..., f_{d-1}$  are algebraically independent.
Then any  element of $ H^0(X, {\mathcal{O}}_X(nD))$
is of the form

$$f=\sum_{i=0}^{m}c_i f_i. 
$$

 Let $P\neq P_0$ be a point on $U$
 and  $x_1,..., x_d$  be the local coordinates at 
 $P$ in an open subset of $V$, $P\in V$.  The conditions that 
 $f(P_0)=f(P)=0$ and $f$ 
 is singular at $P$
   are 
 determined by the 
 system of the following 
linear equations.

$$f(P_0)=c_0f_0(P_0)+c_1f_1(P_0)+...+c_mf_m(P_0)  =0$$
$$f(P)=c_0f_0(P)+c_1f_1(P)+...+c_mf_m(P)  =0$$
  
    $$\frac{\partial f}{\partial x_1}(P)  =c_0\frac{\partial f_0}
    {\partial x_1}(P)+ c_1\frac{\partial f_1}{\partial x_1}(P)+......
    +c_m\frac{\partial f_m}{\partial x_1}(P)=0$$
    
    $$\cdot\cdot\cdot\cdot\cdot\cdot\cdot\cdot\cdot\cdot\cdot\cdot\cdot
\cdot\cdot\cdot\cdot\cdot\cdot\cdot\cdot\cdot
\cdot\cdot\cdot$$ 

  $$\frac{\partial f}{\partial x_d}(P)  =
  c_0\frac{\partial f_0}
    {\partial x_d}(P)+
   c_1\frac{\partial f_1}{\partial x_d}(P)+......
    +c_m\frac{\partial f_m}{\partial x_d}(P) =0.$$

Consider the  $(d+2)$ by $m+1$ matrix $A$

\[\left( \begin{array}{llll}
f_0(P_0) & f_1(P_0) &f_2(P_0)\cdots & f_m(P_0)\\
f_0(P) & f_1(P) &f_2(P)\cdots & f_m(P)\\
 \frac{\partial f_0}
    {\partial x_1}(P)  & \frac{\partial f_1}{\partial x_1}(P)  & \frac{\partial f_2}{\partial x_1}(P)\cdots & \frac{\partial f_m}{\partial x_1}(P)\\
\frac{\partial f_0}
    {\partial x_2}(P)  & \frac{\partial f_1}{\partial x_2}(P) & \frac{\partial f_2}{\partial x_2}(P) \cdots & \frac{\partial f_m}{\partial x_2}(P)\\
\multicolumn{4}{c}\dotfill\\
\frac{\partial f_0}{\partial x_d}(P) & \frac{\partial f_1}{\partial x_d}(P) & \frac{\partial f_2}{\partial x_d} (P) \cdots & \frac{\partial f_m}{\partial x_d}(P)
\end{array}   \right). \]

Since $f_0, ...,f_{d-1}$ are algebraically independent, the determinant of 
the Jacobian 
 $J(\frac{\partial f_i}{\partial x_j})_{0\leq i, j\leq (d-1)}$

\[\left( \begin{array}{lll}
 \frac{\partial f_0}{\partial x_1}  & \frac{\partial f_1}{\partial x_1}\cdots & \frac{\partial f_d}{\partial x_1}\\
\frac{\partial f_0}{\partial x_2} & \frac{\partial f_1}{\partial x_2} \cdots & \frac{\partial f_d}{\partial x_2}\\
\multicolumn{3}{c}\dotfill\\
 \frac{\partial f_0}{\partial x_d} & \frac{\partial f_1}{\partial x_d} \cdots & \frac{\partial f_d}{\partial x_d}
\end{array}   \right) \]

is not identically zero (\cite{R}, Chapter 6, Proposition 6A.4).

Because  $\Phi_{|nD|}$ is a birational  morphism and an 
isomorphism on $U$, the two vectors
$(f_0(P_0), f_1(P_0),..., f_m(P_0))$
and 
 $(f_0(P), f_1(P),..., f_m(P))$
are linearly independent since $P_0\in U$. 
Thus the rank of 
matrix $A$  is at least 2 at any point $P\neq P_0$ of  $X$
except these finite points. At $P_0$, the Jacobian 
$J(\frac{\partial f_i}{\partial x_j})_{1\leq i, j\leq (d-1)}(P_0)$
has rank $d$.

If the rank of $A$ at $P\in X$ is 2, then 
$$\frac{\partial f_i}{\partial x_j}(P)=0, 
\mbox{for\quad  all}\quad  i=0,1,..., d-1;  j=1,..., d. $$

Let $X_2$ be the set on $U$  such that 
the rank of $A$ at every point $P$ of
$U$ is 2. Then the dimension of 
$X_2$ is at most 0. And the dimension of the projective space $C_2$ of  solutions
$(c_0,c_1,..., c_m)$ is  $m-2$ in ${\mathbb{P}^m}$.
So dim$X_2+$dim$C_2=m-2$.

If the rank of $A$ at the point $P$ is 3, then three  rows
including the first two are linearly independent
as vectors in $k^m$. So other $d-1$ rows can be written
as  linear combinations these two rows. There are $d-1$ equations.
Let $X_3$ be the set of points in $U$ such that at every point 
$P$  of $X_3$, rank of $A$ is 3. Then the dimension of 
$X_3$ is  at most 1. Let $C_3$ be the corresponding set of
solutions of the system, then the  dimension of $C_3$ 
as a projective space is
$m-3$. So  dim$C_3+$dim$X_3=m-2$.

In general, if the  rank of $A$  is $r$, $2\leq r\leq d+2$, then 
the rank of the Jacobian $J(f_1,..., f_d)$
is $r-2$, which give $d-(r-2)$ conditions. So
the dimension of $X_r$ is $r-2$. Again we have 
dim$X_r+$dim$C_r=(r-2)+(m-r)=m-2$.

Let $S_x$  be the set  of  effective divisors  $E$
in $|nD|_{P_0}$  such that  $E$  is singular at $x$
and  $S=\{<x,  E>| x\in  X,  E\in  S_x\}$. 
Let  $L=\{f\in H^0(X, {\mathcal{O}}_X(nD))| f(P_0)=0 \}$.
As a projective space, $L$  has dimension 
$m-1$.   
Consider the projections $p_1: S\rightarrow X$
and $p_2: S\rightarrow L$. $p_1$  is surjective.
Take an irreducible 
component $S_1$ of $S$  such that 
$p_1$  is surjective on $S_1$.

There are finitely many  affine  open  subsets  $\{U_i\}$
covering  $X$. 
Since the projective dimension of every fiber over $U_i$
is   $m-d-2$, the dimension of $S_1$ is at most
$m-d-2+d=m-2$. So every irreducible component of $S$  has 
dimension at most $m-2$. 
For the projection $p_2: S\rightarrow L$,
we have dim$p_2(S)\leq $dim$S\leq$ dim$L-1$.
But the projective dimension of $L$ is  $m-1$, 
$p_2(S)$  is a closed subset of $M$. Thus a  general member of
$|nD|_{P_0}$  is smooth.

\begin{flushright}
 Q.E.D. 
\end{flushright}

\section{Applications}

\begin{theorem} Let $M$  be a compact connected  complex manifold
biholomorphic to  an irreducible smooth projective variety $X$
in ${\mathbb{P}}^n({\mathbb{C}})$. Then 

(1) the analytic inverse image of  a general hyperplane
 passing through   a point $P_0$  in ${\mathbb{P}}^n({\mathbb{C}})$
 is a connected complex manifold of codimension 1 on $M$;
 
 (2) the analytic inverse image of  a general hypersurface
 passing through   $q\leq n+1$ points $P_0, ..., P_{q-1}$
 in general position   in ${\mathbb{P}}^n({\mathbb{C}})$
 is a connected complex manifold of codimension 1 on $M$.
\end{theorem}

$Proof$. This is a direct consequence of Theorem 
2.5 and Theorem 3.1.

\begin{flushright}
 Q.E.D. 
\end{flushright}

\begin{theorem} Let $M$  be a compact connected  complex manifold
and $f:  M\rightarrow  X$  is   a  proper holomorphic 
surjective map with maximal rank at every point of $M$,
 where $X$  is a smooth projective variety in
${\mathbb{P}}^n({\mathbb{C}})$. Then 

(1) the analytic inverse image of  a general hyperplane
 passing through   a point $P_0$  in ${\mathbb{P}}^n({\mathbb{C}})$
 is a connected complex manifold of codimension 1 on $M$;
 
 (2) the analytic inverse image of  a general hypersurface
 passing through   $q\leq n+1$ points $P_0, ..., P_{q-1}$
 in general position   in ${\mathbb{P}}^n({\mathbb{C}})$
 is a connected complex manifold of codimension 1 on $M$.
\end{theorem}

Let  $L$  be  a holomorphic line bundle on a 
complex manifold $M$. If for every  $x\in M$,
there is a section $s\in H^0(M, L)$ such that 
$s(x)\neq 0$, then the basis of 
$H^0(M, L)$  gives a holomorphic map  $f$
from $M$
to ${\mathbb{P}}^n({\mathbb{C}})$. If this map 
is an (analytic)  isomorphism  from $M$  to its 
image, then  $L$  is very ample. 

For any  $q\leq n+1$  points $P_0,..., P_{q-1}$
on  $M$, we say that they are in general position if 
their images   under the 
map $f$  given by basis of  $H^0(M, L)$
are in general position in ${\mathbb{P}}^n({\mathbb{C}})$

\begin{corollary} Let $M$  be a compact connected  complex manifold
and  $L$  a very ample holomorphic line bundle 
on  $M$. Then 

(1) for any point $P_0$  on $M$, a general  section $s\in H^0(M, L)$  
with $s(P_0)=0$ gives 
  a connected complex manifold of codimension 1 on $M$;
 
 (2) for any $q\leq  n+1$  points 
 $P_0,..., P_{q-1}$ in general position on
 $M$,  a general  section $s\in H^0(M, L^{\otimes a})$  passing through these 
 $q$  points  gives 
  a connected complex manifold of codimension 1 on $M$, where $a>1$.
\end{corollary}

\begin{theorem}
Let  $Y$  be an irreducible smooth affine variety  in $k^n$
contained in an irreducible  smooth projective variety  $X$
in ${\mathbb{P}}^n$. 
 Let $P_0, P_1,..., P_{q-1}$  be 
$q\leq n+1$  closed points in general position in 
$k^n$. Then for any degree 
$a>1$, a general hypersurface in $k^n$ passing though 
these $q$  points are irreducible and smooth. 

\end{theorem}

$Proof.$ If all $q$ points lie on $Y$, it is a direct 
consequence of Theorem 3.1. In other cases, the proof is similar to
 the proof of  Theorem 3.1.

\begin{flushright}
 Q.E.D. 
\end{flushright}

The following proposition should not be new.
The idea can be traced back to 
Bertini and Severi (\cite{K1}, Section 3 and 4).
I cannot find a proof 
anywhere so I will write one 
here. It is interesting relationship 
between the linear system of a divisor and 
linear  system of 
hypersurfaces.

\begin{proposition} Let $D$  be an effective Cartier  divisor
on an irreducible   projective variety $X$ in ${\mathbb{P}}^N$ such
that  $\kappa(D, X)\geq 1$. Then there is a positive
integer $n_0$ such that for all $n\geq n_0$, 
the linear system $|nD|$, except for the divisor 
$nD$,  can be obtained by cutting out on $X$
by a linear system of hypersurfaces and then removing some fixed 
components, which are the common components of all hypersurfaces
in the system.   
\end{proposition}

$Proof$. Since $s=\kappa(D, X)\geq 1$, the dimension of 
the linear system $|nD|$ as a vector space grows like $cn^s$,
where $c>0$ is a constant. Let  $\{f_0, f_1,..., f_m\}$  be a basis 
of 
$H^0(X, {\mathcal{O}}_X(nD))$, then  there are rational functions $g_i$ on $X$
such that div$(f_i)=$div$(g_i)+nD$.  Let $g_i=h_i/h_i'$, where  
$h_i$  and $h_i'$ are homogeneous polynomials of the same degree
in ${\mathbb{P}}^N$. Any element $g$ of 
$L$ is a linear combination of 
its basis. 
Let  ${\mathbb{C}}(X)$ 
be the function field of $X$, then 
the vector space 
$H^0(X, {\mathcal{O}}_X(nD))$
is isomorphic to the vector space  (up to a constant)
$$L=\{g\in {\mathbb{C}}(X)|g=0\quad  {\mbox{or}}\quad {\mbox{div}}(g)+nD\geq 0\}.
$$
Since $|nD|$ is isomorphic to 
$L$ module a constant
(\cite{U}, Chapter II, Lemma 4.16),  
 there are constants 
$c_0, c_1,..., c_n$  such that
 we can write  any element $E\in |nD|$ as follows

$$   E={\mbox{div}} (\sum_{i=0}^m c_i g_i)+nD
\quad\quad\quad\quad\quad\quad\quad\quad\quad\quad \quad\quad\quad\quad\quad \quad\quad\quad\quad\quad
\quad\quad\quad$$

$$
={\mbox{div}}(\sum_{i=0}^m c_i  (\frac{h_i}{h_i'}))+nD
\quad\quad\quad\quad\quad\quad\quad\quad\quad\quad \quad\quad\quad\quad\quad \quad\quad\quad\quad\quad\quad$$ 

 $$
={\mbox{div}}(\sum_{i=0}^m c_i (\frac{h_0'...h_{i-1}'h_ih_{i+1}'...h_m'}
{h_0'...h_{i-1}'h_i'h_{i+1}'...h_m'}))+nD\quad\quad\quad\quad\quad \quad\quad\quad\quad\quad\quad\quad\quad
$$

$$
={\mbox{div}} (\sum_{i=0}^m c_i  (h_0'...h_{i-1}'h_ih_{i+1}'...h_m'))
-\mbox{div}(h_0'...h_{i-1}'h_i'h_{i+1}'...h_m')
+nD\quad
$$

 $$
={\mbox{div}} ( \sum_{i=0}^m c_i  (\alpha_i))
-\mbox{div}(\beta)+
nD,\quad\quad\quad\quad\quad\quad\quad\quad\quad\quad \quad\quad\quad\quad\quad \quad
$$

where 
$$\alpha_i=h_0'...h_{i-1}'h_ih_{i+1}'...h_m'$$
 and 
$$\beta=h_0'...h_{i-1}'h_i'h_{i+1}'...h_m'$$
 are homogeneous polynomials in th projective space 
 ${\mathbb{P}}^N$.
 From the above formula we see that $\beta$ defines a fixed 
 hypersurface.

The equation  proves the proposition.

\begin{flushright}
 Q.E.D. 
\end{flushright}

\begin{center} {\bf Acknowledgments}

\end{center}

I thank Professor Robin Hartshorne who  suggested me to consider
multiple points and pointed out to me that the points may not 
lie on the variety. I have  benefited from  the communication 
with  Professor Steven L. Kleiman. 
I am very grateful for his  time  and patience.
 I also thank   Professor Igor Dolgachev
who allows me to include Theorem 2.5.

 \end{document}